\documentclass{amsart}
\newcommand\str{{\text{\raisebox{0.5ex}{\normalfont\bf\Large.}}}}
\usepackage{latexsym}
\usepackage[all]{xypic}
\usepackage{isolatin1}
\numberwithin{equation}{section}
\newtheorem{Lem}{Lemma}[section]
\newtheorem{Prop}[Lem]{Proposition}
\newtheorem{Cor}[Lem]{Corollary}
\newtheorem{Thm}[Lem]{Theorem}
\theoremstyle{definition}
\theoremstyle{remark}
\newtheorem{Rem}[Lem]{Remark}
\newcommand{\mapsfrom}{\mathrel{\mbox{$\leftarrow\joinrel\mapstochar$\,}}}
\renewcommand\o{\otimes}
\DeclareMathOperator\Hom{\operatorname{Hom}}

\DeclareMathOperator\ev{\operatorname{ev}}
\DeclareMathOperator\db{\operatorname{db}}

\newcommand\bop{{\operatorname{bop}}}
\newcommand\drdot[1]{{^\str_\str #1^\str}}
\newcommand\HMod[4]{{^{#1}_{#3}\mathcal M^{#2}_{#4}}}
\newcommand\LMod[1]{{_{#1}\mathcal M}}
\newcommand\LComod[1]{{^{#1}\mathcal M}}
\newcommand\RComod[1]{\mathcal M^{#1}}
\newcommand\RComodf[1]{\mathcal M^{#1}_f}
\newcommand\RMod[1]{\mathcal M_{#1}}
\newcommand\LModf[1]{{_{#1}\mathcal M_f}}
\newcommand\BiMod[1]{{_{#1}\mathcal M_{#1}}}
\newcommand\BiComod[1]{{^{#1}\mathcal M^{#1}}}
\newcommand\kmod{\mathcal M_k}
\newcommand{\ou}[1]{\mathrel{\mathop{\otimes}_{#1}}}
\newcommand{\co}[1]{\mathrel{\mathop{\Box}_{#1}}}

\newcommand\si[1]{(#1)}
\newcommand\sw[1]{{}_{(#1)}}
\newcommand\swm[1]{{}_{(-#1)}}
\newcommand\so[1]{^{(#1)}}
\newcommand\som[1]{^{(-#1)}}
\newcommand\ol{\overline}
\newcommand\inv{^{-1}}
\renewcommand\epsilon\varepsilon

\newcommand\xihut{\hat\xi}
\newcommand\pEv{E'}
\newcommand\Ev{E}
\newcommand\Db{D}
\newcommand\tgamma{\widetilde\gamma}
\newcommand\ttheta{\widetilde\vartheta}
\newcommand\tttheta{\vartheta^t}
\newcommand\imp{\ensuremath{\implies}}
\makeatletter
\def\namelabel#1#2{\@bsphack
  \protected@write\@auxout{}%
         {\string\newlabel{#1.nme}{{#2}{#2}}}%
  \@esphack}
\def\nmlabel#1#2{\label{#2}\namelabel{#2}{#1}}
\newcommand\nmref[1]{\ref{#1.nme}\ \ref{#1}}
\makeatother
\begin{document}
\title{Two characterizations of finite quasi-Hopf algebras}
\author{Peter Schauenburg}
\address{Mathematisches Institut der Universit\"at M\"unchen, 
Theresienstr.~39, 80333~M\"unchen, Germany}
\email{schauen@rz.mathematik.uni-muenchen.de}
\subjclass{16W30}
\keywords{Quasi-Hopf algebra}
\begin{abstract}
  Let $H$ be a finite-dimensional quasibialgebra. We show that 
  $H$ is a quasi-Hopf algebra if and only if the monoidal 
  category of its finite-dimensional left modules is rigid, 
  if and only if a structure theorem for Hopf modules
  over $H$ holds. We also show that a dual structure theorem for
  Hopf modules over a coquasibialgebra $H$ holds if and only if
  the category of finite-dimensional right $H$-comodules is rigid;
  this is not equivalent to $H$ being a coquasi-Hopf algebra.
\end{abstract}
\maketitle
\section{Introduction}
Let $H$ be a bialgebra over the field $k$. 
It was shown by Ulbrich \cite{Ulb:HARMC} that $H$ is a Hopf algebra
if and only if the monoidal category $\RComodf H$ of finite-dimensional
right $H$-comodules is rigid, that is, if every finite-dimensional
$H$-comodule has a dual object within the category $\RComodf H$.
In particular, if $H$ is a finite-dimensional bialgebra, then 
$H$ is a Hopf algebra if and only if the category $\LModf H$
of finite-dimensional left $H$-modules is rigid.

It is a natural question whether the same holds for quasibialgebras:
Drinfeld's definition of a quasibialgebra $H$ ensures that the category
$\LMod H$ is, just like in the bialgebra case, a monoidal category.
And Drinfeld's definition of a quasiantipode is motivated by the fact
that the category of finite-dimensional modules over a quasi-Hopf
algebra is a rigid monoidal category. However,
if we try to prove the converse, then 
we run into difficulties. The key problem is that the underlying
functor $\LMod H\rightarrow\kmod$ to the category of $k$-vector
spaces is monoidal if $H$ is a bialgebra, and monoidal functors
automatically preserve dual objects. If $H$ is only a quasibialgebra,
then the underlying functor $\LMod H\rightarrow \kmod$ is still
compatible with tensor products, but not coherent in the sense of
the definition of a monoidal functor. Thus it is not clear that
the functor preserves duals. That the problem is really
serious was shown in \cite{Sch:HAEMC} by an example based on a
construction of Yongchang Zhu \cite{Zhu:HARRHA}: There is a coquasibialgebra $H$
such that the category $\RComodf H$ is rigid, although $H$ is
not a coquasi-Hopf algebra. The existence of a coquasiantipode
is ruled out quite drastically by the fact that a finite-dimensional
$H$-comodule and its dual object may have different dimensions.
One result of this paper will be that all is well in the 
finite-dimensional case: A finite-dimensional quasibialgebra $H$ is a 
quasi-Hopf algebra if and only if $\LModf H$ is rigid.

Another well-known criterion says that a bialgebra $H$ is 
a Hopf algebra if and only if the structure theorem for Hopf modules
holds, that is, if the obvious functor
$\kmod\rightarrow\HMod{}H{}H$ mapping a vector space $V$ to
$V\o H^\str_\str$ is a category equivalence.
If we try to establish a version of this criterion for quasibialgebras,
the first problem is that there are no Hopf modules: If $H$ is a
quasibialgebra, then it is not a coassociative coalgebra, so one does
not know what a comodule should be. This first problem was solved by
Hausser and Nill \cite{HauNil:ITQHA}, who observed that one
can still define a Hopf (bi)module category $\HMod{}HHH$; we may
say briefly that it is the category of $H$-comodules over the
coassociative coalgebra $H$ {\em within the monoidal category 
$\BiMod H$.}
Moreover, Hausser and Nill prove a structure theorem for Hopf
bimodules: A certain functor $\LMod H\rightarrow\HMod{}HHH$
is a category equivalence if $H$ is a quasi-Hopf algebra. If we
try to prove a converse, we run into difficulties once again.
In the case of ordinary bialgebras, the proof is based on 
another criterion: $H$ is a Hopf algebra if the canonical 
map $H\o H\ni g\o h\mapsto gh\sw 1\o h\sw 2\in H\o H$ is a bijection.
If the structure theorem for Hopf modules holds, it is very easy to
check that the canonical map is bijective. For a quasi-Hopf algebra
$H$, Drinfeld's paper \cite{Dri:QHA} contains an analog
of the canonical map $H\o H\rightarrow H\o H$. However,
the proper anolog is given by a more complicated formula;
in particular,  one already
needs a quasiantipode to even write down the map (or its inverse),
and it seems to have no analog for quasibialgebras. 
We shall show that the problem is just as serious as
that with the first criterion mentioned above: 
Given a coquasibialgebra $H$ we will prove
in \nmref{sec:duho} that the structure theorem
for Hopf modules holds --- that is, a certain functor
$\RComod H\rightarrow\HMod HHH{}$ is a a category equivalence ---, if
and only if the category $\RComodf H$ is rigid. 
For a finite-dimensional quasi-Hopf algebra this provides
a new and rather conceptual proof for the structure theorem of
Hausser and Nill --- in fact our proof needs hardly any 
unpleasant calculations
with the quasibialgebra structure and its axioms, and none at all
with the quasiantipode. 
Of course the result also provides examples of coquasibialgebras
that satisfy the structure theorem for Hopf modules, while they
do not have a coquasiantipode.

On the other hand, if $H$ is a finite-dimensional quasibialgebra,
\nmref{char} shows that all is well, that is, the structure theorem
for Hopf modules is equivalent to the existence of a quasiantipode.

In \nmref{sec:bij} we will show that the quasiantipode
of a finite-dimensional quasi-Hopf algebra is a bijection. This
was first proved by Bulacu and Caenepeel \cite{BulCae:IDQHAA}. 
We will give a rather different proof.

Throughout the paper, we work over a base field $k$.

\section{Duality and the structure of Hopf modules}\nmlabel{Section}{sec:duho}

Throughout this section, we let $H$ be a coquasibialgebra.
That is, $(H,\Delta,\epsilon)$ is a coassociative coalgebra,
endowed with a (nonassociative) multiplication $\nabla\colon H\o H\rightarrow H$
which is a coalgebra map, a grouplike unit element $1\in H$, and
a convolution invertible form $\phi\colon H\o H\o H\rightarrow k$,
the associator, satisfying the identities
$\phi(g\o 1\o h)=\epsilon(g)\epsilon(h)$,
$$  (f\sw 1g\sw 1)h\sw 1\phi(f\sw 2\o g\sw 2\o h\sw 2)
    =\phi(f\sw 1\o g\sw 1\o h\sw 1)f\sw 2(g\sw 2h\sw 2),$$
and
\begin{multline*}
  \phi(d\sw 1f\sw 1\o g\sw 1\o h\sw 1)\phi(d\sw 2\o f\sw 2\o g\sw 2h\sw 2)
    \\=\phi(d\sw 1\o f\sw 1\o g\sw 1)\phi(d\sw 2\o f\sw 2g\sw 2\o h\sw 1)\phi(f\sw 3\o g\sw 3\o h\sw 2)
\end{multline*}
for $d,f,g,h\in H$. We have used Sweedler notation in the form 
$\Delta(h)=h\sw 1\o h\sw 2$. We will also use Sweedler notations
$V\ni v\mapsto v\sw 0\o v\sw 1\in V\o H$ for right, and 
$V\ni v\mapsto v\swm 1\o v\sw 0\in H\o V$ for left comodule structures.

The axioms ensure that the category $\RComod H$ of right $H$-comodules
is a monoidal category in the following way: For $V,W\in\RComod H$,
the tensor product $V\o W$ over $k$ is an $H$-comodule with the
codiagonal comodule structure induced by multiplication in $H$;
the associativity isomorphism 
$\Phi\colon (U\o V)\o W\rightarrow U\o (V\o W)$ is given
by $\Phi(u\o v\o w)=u\sw 0\o v\sw 0\o w\sw 0\phi(u\sw 1\o v\sw 1\o w\sw 1)$.
Since the opposite of a coquasibialgebra and the tensor product
of two coquasibialgebras are coquasibialgebras as well, the category
$\BiComod H$ of $H$-$H$-bicomodules is also a monoidal category.
This time the associativity isomorphism 
$\Phi\colon (U\o V)\o W\rightarrow U\o(V\o W)$ is
given by
$\Phi(u\o v\o w)=\phi\inv(u\swm 1\o v\swm 1\o w\swm 1)u\sw 0\o v\sw 0\o w\sw 0\phi(u\sw 1\o v\sw 1\o w\sw 1)$.
It is a key observation that $H$ (which is not associative as
a $k$-algebra) is an associative algebra within the monoidal
category $\BiComod H$, that is, we have
$$\nabla(\nabla\o H)=\nabla(H\o\nabla)\Phi\colon (H\o H)\o H\rightarrow H\o (H\o H).$$
Thus we can use the general theory of algebras and modules
in monoidal categories, see Pareigis \cite{Par:NARMTI,Par:NARMTII}, 
to do (or rather avoid) calculations with the multiplicative
structure of $H$.

In particular, there is a well-defined notion of (say, left)
$H$-module within the monoidal category $\BiComod H$. We denote
the category of such modules by $\HMod HHH{}$, and call its objects
Hopf modules. We note that for any $M\in\HMod HHH{}$ and $P\in\BiComod H$
we have $M\o P\in\HMod HHH{}$ with the ``obvious''
left module structure
$$H\o(M\o P)\xrightarrow{\Phi\inv}(H\o M)\o P\xrightarrow{\mu\o P}M\o P,$$
where $\mu$ denotes the module structure of $M$.
We will abbreviate this module structure by a dot attached to 
the tensorand $M$, i.e. write $\drdot M\o{^\str P^\str}$ for it,
with the upper dots indicating on which tensor factors we have
a codiagonal coaction, and the upper dot indicating where the 
action takes place; note, though, that the actual formula for
the action involves both tensorands through the action of 
the associators.

Taking $M=H$ as a special case we obtain the left adjoint
$P\mapsto \drdot H\o{^\str P^\str}$ to the underlying 
functor $\HMod HHH{}\rightarrow\BiComod H$.

As a particular case, we can consider a right $H$-comodule $V$
as a bicomodule with the trivial comodule structure on the left,
and apply the above construction to obtain a functor
$$\mathcal L\colon\RComod H\ni V\mapsto {^\str_\str H^\str}\o V^\str\in\HMod HHH{}.$$
The formally dual version of this functor (for a quasibialgebra)
was studied by Hausser and Nill \cite{HauNil:ITQHA}, who also 
proved that it is an equivalence in case $H$ is a quasi-Hopf algbebra.
Moreover, $\HMod HHH{}$ is a monoidal category, and
$\mathcal L$ is a monoidal functor. Hausser and Nill show this in 
the dual case using the assumption that $H$ is quasi-Hopf, the
quasibialgebra case is treated in \cite{Sch:HMDQHA}. We shall say
for short that the structure theorem for Hopf modules holds if
$\mathcal L$ is a category equivalence. In this section we shall give
a different proof of the structure theorem for Hopf modules
than the ones in \cite{HauNil:ITQHA,Sch:HMDQHA}, under the weaker
assumption that $\RComodf H$ is a rigid monoidal category.

We start by stating several facts on the functor $\mathcal L$
and the category $\HMod HHH{}$ that are formally dual to 
facts proved and used in \cite{HauNil:ITQHA} and \cite{Sch:HMDQHA}. 
We will not give the
formally dual proofs, but will indicate how
\nmref{monfun} follows from more abstract reasons without
any work.
\begin{Lem}[dual to part of {\cite[Prop.3.6]{Sch:HMDQHA}}]
  The functor $\mathcal L\colon\RComod H\rightarrow \HMod HHH{}$
  is exact, fully faithful, and commutes with arbitrary colimits.
  In particular, colimits
  and equalizers of diagrams whose objects are in the image of
  $\mathcal L$ are also in that image.
\end{Lem}
The dual statement of the following is observed between Corollary
3.9 and Lemma 3.10 of \cite{HauNil:ITQHA}. See also
\cite[Lem.and Def.3.2]{Sch:HMDQHA}.
\begin{Lem}
  The category $\HMod HHH{}$ is a monoidal category in 
  the following way:

  The tensor product of $M,N\in\HMod HHH{}$ is their cotensor 
  product $M\co HN$ equipped with the module structure
  given by $h(m\o n)=h\sw 1m\o h\sw 2n$.

  In particular the underlying functor $\HMod HHH{}\rightarrow\BiComod H$
  is a strict monoidal functor with the monoidal category structure
  on the target given by cotensor product.
\end{Lem}
Dually to \cite[Lem.3.4]{Sch:HMDQHA} one can check that
for any $M\in \HMod HHH{}$ and $V\in\RComod H$ 
the canonical
isomorphism
\begin{align*}
  \xihut\colon M\co H(H\o V)&\rightarrow M\o V\\
  m\sw 0\o m\sw 1\o v&\mapsto m\o v\\
  m\o h\o v&\mapsfrom m\epsilon(h)\o v
\end{align*}
is a morphism in $\HMod HHH{}$. This follows from the following
more general statement:
\begin{Lem}
  Let $M,N\in\HMod HHH{}$, and $V\in\RComod H$. The canonical 
  isomorphism (``the identity'')
  $$(M\co HN)\o V\cong M\co H(N\o V)$$
  is an isomorphism in $\HMod HHH{}$. If we identify
  $(M\co HN)\o V=M\co H(N\o V)=M\co HN\o V$, then 
  $$\Phi_{M\co HN,V,W}=M\co H\Phi_{N,V,W}\colon (M\co HN\o V)\o W\rightarrow M\co H N\o(V\o W).$$
\end{Lem}
\begin{proof}
  It is obvious that the isomorphism is left and right $H$-colinear.
  $H$-linearity is a small calculation: Denoting the respective
  actions by $h((m\o n)\o v)$ and $h(m\o(n\o v))$ for $m\o n\in M\co HN$,
  $v\in V$, and $h\in H$, we find
  \begin{multline*}
    h((m\o n)\o v)=h\sw 1(m\o n)\sw 0\o v\sw 0\phi(h\sw 2\o(m\o n)\sw 1\o v\sw 1)
      \\=h\sw 1(m\o n\sw 0)\o v\sw 0\phi(h\sw 2\o n\sw 1\o v\sw 1)
      \\=h\sw 1m\o h\sw 2n\sw 0\o v\sw 0\phi(h\sw 3\o n\sw 1\o v\sw 1)
      \\=h\sw 1m\o h\sw 2(n\o v)
      =h(m\o(n\o v))
  \end{multline*}
The two associativity isomorphisms both map
$m\o n\o v\o w$ to $m\o n\sw 0\o v\sw 0\o w\sw 0\phi(n\sw 1\o v\sw 1\o w\sw 1)$.
\end{proof}
We can restate the Lemma as follows: The category $\HMod HHH{}$ 
of left $H$-modules in $\BiComod H$ is naturally a right $\BiComod H$-category,
hence a right $\RComod H$-category. On the other hand $\HMod HHH{}$
is a monoidal category, so it is naturally a left $\HMod HHH{}$-category.
The Lemma says that the right $\RComod H$-category structure is 
compatible with the left $\HMod HHH{}$-category structure (or we
have an $\HMod HHH{}$-$\RComod H$-bicategory).
It follows as in \cite[Thm.3.3]{Sch:AMCGHSP} that a monoidal 
functor $(\mathcal L,\xi)\colon\RComod H\rightarrow\HMod HHH{}$
is given by
$\mathcal L(V)=H\o V$, and 
$$\xi\colon (H\o V)\co H(H\o W)\xrightarrow{\xihut} (H\o V)\o W\xrightarrow\Phi H\o (V\o W).$$
We shall repeat briefly the abstract argument: The right
action of $\RComod H$ on $\HMod HHH{}$ gives a monoidal functor
(here actually an antimonoidal functor) from $\RComod H$
to the monoidal category of endofunctors of $\HMod HHH{}$.
Since the action is compatible with the action of $\HMod HHH{}$, 
the antimonoidal functor has its image in the category of
$\HMod HHH{}$-endofunctors of $\HMod HHH{}$, which is (anti)monoidally
equivalent to the monoidal category $\HMod HHH{}$ itself (compare
to the fact that the endomorphism ring of a ring $R$ considered
as an $R$-module, is isomorphic to the ring itself). Thus we have
a monoidal functor fron $\RComod H$ to $\HMod HHH{}$.

While this abstract argument is considerably simpler than the 
calculations needed in either \cite{HauNil:ITQHA} or
\cite{Sch:HMDQHA} to show that $\mathcal L$ is monoidal, 
the reader might not want to get into $\mathcal C$-category theory.
In this case one may dualize the proof of \cite[Prop.3.6]{Sch:HMDQHA}
to obtain:
\begin{Prop}\nmlabel{Proposition}{monfun}
  Define
  $$\xi=\xi_{VW}=\left((H\o V)\co H(H\o W)\xrightarrow{\xihut}(H\o V)\o W\xrightarrow\Phi H\o(V\o W)\right)$$
  for $V,W\in \RComod H$. Then
  $(\mathcal L,\xi)\colon\RComod H\rightarrow\HMod HHH{}$ is a monoidal
  functor.
\end{Prop}

Let $\mathcal C$ be a monoidal category. Recall that a dual object
of $V\in\mathcal C$ is a triple $(V^\vee,\ev,\db)$ in which 
$V^\vee\in\mathcal C$, and $\ev\colon V^\vee\o V\rightarrow I$
and $\db\colon I\rightarrow V\o V^\vee$ are morphisms such that
the two compositions
\begin{gather*}
  V\xrightarrow{\db\o V}(V\o V^\vee)\o V\xrightarrow\Phi V\o(V^\vee\o V)\xrightarrow{V\o\ev}V\\
  V^\vee\xrightarrow{V^\vee\o\db}V^\vee\o(V\o V^\vee)\xrightarrow{\Phi\inv}(V^\vee\o V)\o V^\vee\xrightarrow{\ev\o V^\vee}V^\vee
\end{gather*}
are identities. If $(\mathcal F,\xi)\colon \mathcal C\rightarrow \mathcal D$
is a monoidal functor, and $(V^\vee,\ev,\db)$ is a dual object
of $V$ in $\mathcal C$, then 
$\mathcal F(V^\vee)$ is a dual object of $\mathcal F(V)$ in $\mathcal D$,
with evaluation and coevaluation
\begin{gather*}
  \mathcal F(V^\vee)\o\mathcal F(V)\xrightarrow\xi\mathcal F(V^\vee\o V)\xrightarrow{\mathcal F(\ev)}\mathcal F(I)\cong I\\
  I\cong\mathcal F(I)\xrightarrow{\mathcal F(\db)}\mathcal F(V\o V^\vee)\xrightarrow{\xi\inv}\mathcal F(V)\o\mathcal F(V^\vee).
\end{gather*}

Let $V\in\RComod H$ be finite-dimensional. We can endow the dual
vector space $V^*$ with a canonically corresponding left $H$-comodule
structure defined by
$\varphi\swm 1\varphi\sw 0(v)=\varphi(v\sw 0)v\sw 1$
for all $\varphi\in V^*$ and $v\in V$. Equivalently,
$v_{i\si 0}\o v_{i\si 1}\o v^i=v_i\o v^i\swm 1\o v^i\sw 0\in V\o H\o V^*$,
that is $v_i\o v^i\in V\co H V^*$. Note that the map
$$\pEv\colon{^\str(V^*)}\o V^\str\ni \varphi\o v\mapsto \varphi(v\sw 0)v\sw 1=\varphi\swm 1\varphi\sw 0(v)\in H$$
is an $H$-bimodule map.
\begin{Lem}\nmlabel{Lemma}{duallemma}
  Let $V\in\RComodf H$. Then a dual object of $\mathcal L(V)$
  in the monoidal category $\HMod HHH{}$ is given by
  $$({^\str_\str H^\str}\o {^\str(V^*)},\Ev,\Db)$$
  with 
  \begin{gather*}
  \Ev=\left((\drdot H\o {^\str(V^*)})\co H(\drdot H\o {^\str V})\xrightarrow{\xihut}H\o V^*\o V\xrightarrow{H\o \pEv}H\o H\xrightarrow\nabla H\right)\\
  \Db\colon H\ni h\mapsto h\sw 1\o v_i\o h\sw 2\o v^i\in (\drdot H\o {V^\str})\co H(\drdot H\o {^\str(V^*)})
  \end{gather*}
\end{Lem}
\begin{proof}
From the definition of $\Ev$ it is clear that $\Ev$ is 
a well-defined morphism in $\HMod HHH{}$. Note that we have
$\xihut\inv(h\o\varphi\o v)=h\sw 1\o\varphi\o h\sw 2\o v$, and
$\Ev(h\sw 1\o\varphi\o h\sw 2\o v)=h\varphi(v\sw 0)v\sw 1=h\varphi\swm 1\varphi\sw 0(v),$
hence $\epsilon\Ev(h\sw 1\o\varphi\o h\sw 2\o v)=\epsilon(h)\varphi(v)$.
The map $\Db$ is well-defined since $v_i\o v^i\in V\co HV^*$
and $h\sw 1\o h\sw 2\in H\co HH$. It is obviously left and right
$H$-colinear, and it is $H$-linear by the calculation
\begin{multline*}
  g\Db(h)
    =g(h\sw 1\o v_i\o h\sw 2\o v^i)
    =g\sw 1(h\sw 1\o v_i)\o g\sw 2(h\sw 2\o v^i)
    \\=g\sw 1h\sw 1\o v_{i\si 0}\phi(g\sw 2\o h\sw 2\o v_{i\si 1})\o g\sw 3(h\sw 3\o v^i)
    \\=g\sw 1h\sw 1\o v_i\o\phi(g\sw 2\o h\sw 2\o v^i\swm 1)g\sw 2(h\sw 2\o v^i\sw 0)
    \\=g\sw 1h\sw 1\o v_i\o g\sw 2h\sw 2\o v^i
    =\Db(gh)
\end{multline*}
for $g,h\in H$. To check the identities for a dual object, we have
to bear in mind the canonical identifications
\begin{align*}
H\co HM&\cong M&M\co HH&\cong M\\
\sum h_i\o m_i&\mapsto \sum\epsilon(h_i)m_i&\sum m_i\o h_i&\mapsto \sum m_i\epsilon(h_i)\\
m\swm 1\o m\sw 0&\mapsfrom m&m\sw 0\o m\sw 1&\mapsfrom m,
\end{align*}
and can calculate
\begin{multline*}
  ((H\o V)\co H\Ev)(\Db\co H(H\o V))(h\o v)
    =((H\o V)\co H\Ev)(\Db(h\sw 1)\o (h\sw 2\o v))
    \\=h\sw 1\o v_i\epsilon\Ev(h\sw 2\o v^i\o h\sw 3\o v)
    =h\sw 1\o v_i\epsilon(h\sw 2v^i(v\sw 0)v\sw 1)
    \\=h\o v_iv^i(v)=h\o v
\end{multline*}
for $h\in H$ and $v\in V$, and
\begin{multline*}
  (\Ev\co H(H\o V^*))((H\o V^*)\co H\Db)(h\o\varphi)
    =(\Ev\co H(H\o V^*)(h\sw 1\o \varphi\o\Db(h\sw 2))
    \\=\epsilon\Ev(h\sw 1\o\varphi\o h\sw 2\o v_i)(h\sw 3\o v^i)
    =h\o \varphi(v_i)v^i=h\o\varphi
\end{multline*}
for $h\in H$ and $\varphi\in V^*$.
\end{proof}
\begin{Thm}\nmlabel{Theorem}{dualthm}
  Let $H$ be a coquasibialgebra. The following are equivalent:
  \begin{enumerate}
    \item The functor $\mathcal L\colon\RComod H\rightarrow\HMod HHH{}$
      is an equivalence.
    \item The category $\RComodf H$ is rigid.
  \end{enumerate}
\end{Thm}
\begin{proof}
  Assume (1), and let $V\in\RComodf H$. Then $\mathcal L(V)$ 
  has a left dual object in $\HMod HHH{}$ by \nmref{duallemma}.
  Since $\mathcal L$ is an equivalence, $V$ has a left dual object
  $V^\vee$ in $\RComod H$. But $V^\vee$ is necessarily finite-dimensional.
  For let $\db\colon k\rightarrow V\o V^\vee$ be the relevant
  coevaluation. Then $\db(1)=\sum_{i=1}^r x_i\o y_i$ for some
  $x_i\in V$ and $y_i\in V^\vee$. The latter generate a finite-dimensional
  subcomodule $U\subset V^\vee$, and it is straightforward to check
  that the map $\db'\colon k\rightarrow V\o U$ induced by $\db$,
  and the restriction $\ev'\colon U\o V\rightarrow k$ of the 
  evaluation $\ev\colon V^\vee\o V\rightarrow k$ make $U$ a dual object
  of $V$, whence $U=V^\vee$.

  Now assume (2). We need to show that $\mathcal L$ is essentially
  surjective. Let $M\in\HMod HHH{}$. Then, calculating within the
  monoidal category $\BiComod H$, we have
  $M\cong H\ou HM$, that is, we have a coequalizer
  $$\xymatrix{H\o (H\o M)\ar@<.5ex>[r]\ar@<-.5ex>[r]
      &H\o M\ar[r]&M}$$
  in $\HMod HHH{}$, in which the first two objects have the form
  $\drdot H\o {^\str P^\str}$ for $P\in\BiComod H$. Since 
  the image of $\mathcal L$ is closed under coequalizers, it
  suffices to check that objects of this form are in the image
  of $\mathcal L$. Now for $P\in\BiComod H$ we have 
  $P\cong P\co HH$, that is, we have an equalizer
  $$\xymatrix{P\ar[r]&{^\str P}\o H^\str\ar@<.5ex>[r]\ar@<-.5ex>[r]&^\str P\o H\o H^\str}$$
  in $\BiComod H$. Since the image of $\mathcal L$ is closed under
  equalizers, it suffices to check that objects of the form
  $\drdot H\o {^\str W}\o{H^\str}$ for $W\in\LComod H$
  are in the image of $\mathcal L$. Since
  $$\drdot H\o{^\str W}\o H^\str\cong  (\drdot H\o{^\str W})\co H(\drdot H\o H^\str),$$
  and the image of $\mathcal L$ is closed under cotensor product, it
  suffices to verify that $\drdot H\o{^\str W}$ is in the image
  of $\mathcal L$, and since $\mathcal L$ preserves colimits, 
  we may assume that $W$ is finite-dimensional. But then we have
  $W\cong V^*$ for some $V\in\RComodf H$, and 
  $\drdot H\o{^\str W}$ is the left dual of $\mathcal L(V)$ in 
  $\HMod HHH{}$. Since monoidal functors preserve duals, it follows
  that $\drdot H\o{^\str W}\cong\mathcal L(V^\vee)$, where 
  $V^\vee$ is a left dual of $V$ in $\RComodf H$.
\end{proof}

Condition (2) in \nmref{dualthm} is fulfilled if $H$ is a 
coquasi-Hopf algebra. We will not recall the axioms of a 
coquasiantipode in any detail, but shall merely say that it involves
an anti-coalgebra endomorphism $S$ of $H$ 
that allows us to endow $V^\vee:=V^*$, the $k$-linear dual of
$V\in\RComodf H$, with an $H$-comodule structure, and extra structure
elements that make $V^\vee$ into a left dual of $V$ in $\RComodf H$.
If $W\cong V^*$ as at the end of the proof of \nmref{dualthm}, 
then we see that $V^\vee\cong W^S$, the right $H$-comodule 
obtained from the left $H$-comodule $W$ along $S$. Thus we have:
\begin{Cor}\nmlabel{Corollary}{cqhacor}
  If $H$ is a coquasi-Hopf algebra, then the functor
  $\mathcal L\colon \RComod H\rightarrow\HMod HHH{}$ is an equvialence.

  If $W\in\LComod H$, then 
  $\drdot H\o{^\str}W\cong\drdot H\o (W^S)^\str\in\HMod HHH{}$.
\end{Cor}
However, there are examples of coquasibialgebras $H$ such that
$\RComodf H$ is right and left rigid, while $H$ does not have
a coquasiantipode; see \cite[Sec.4.5]{Sch:HAEMC}.
\section{Finite Quasi-Hopf algebras}
Recall that a quasibialgebra $H$ is an associative algebra
with an algebra map $\Delta\colon H\rightarrow H\o H$
called comultiplication, an algebra map $\epsilon\colon H\rightarrow k$
that is a counit for $\Delta$, and an invertible element
$\phi\in H\o H\o H$ such that 
\begin{gather}
(\epsilon\o H)\Delta(h)=h=(H\o\epsilon)\Delta(h),\\
(H\o\Delta)\Delta(h)\cdot\phi=\phi\cdot(\Delta\o H)\Delta(h)\label{quasicoass},\\
(H\o H\o\Delta)(\phi)\cdot(\Delta\o H\o H)(\phi)=(1\o\phi)\cdot(H\o\Delta\o H)(\phi)\cdot(\phi\o 1)\label{cocycle},\\
(H\o\epsilon\o H)(\phi)=1
\end{gather}
hold for all $h\in H$.
We will write $\Delta(h)=:h\sw 1\o h\sw 2$, 
$\phi=\phi\so 1\o\phi\so 2\o\phi\so 3$, and
$\phi\inv=\phi\som 1\o\phi\som 2\o\phi\som 3$. 

A finite-dimensional quasibialgebra is the same as the dual of
a finite-dimensional coquasibialgebra (historically, vice versa
would be more to the point).

A quasiantipode $(S,\alpha,\beta)$ for a quasibialgebra $H$ consists
of an anti-algebra endomorphism $S$ of $H$, 
and elements $\alpha,\beta\in H$,
such that 
\begin{align*}
  S(h\sw 1)\alpha h\sw 2&=\epsilon(h)\alpha, & h\sw 1\beta S(h\sw 2)&=\epsilon(h)\beta,\\
  \phi\so 1\beta S(\phi\so 2)\alpha\phi\so 3&=1,&S(\phi\som 1)\alpha\phi\som 2\beta\phi\som 3&=1
\end{align*}
hold in $H$, for $h\in H$. A quasi-Hopf algebra is a quasibialgebra
with a quasi-antipode. Note that we disagree in this definition with
Drinfeld who requires $S$ to be a bijection. We will return to this
in \nmref{sec:bij} where we give another proof for a recent
result of Bulacu and Caenepeel, which says that the antipode of a 
finite-dimensional quasi-Hopf algebra is automatically bijective.

The main result of this section characterizes finite quasi-Hopf algebras
via rigidity of their module category, or the structure theorem for 
Hopf modules. The functor $\mathcal R$ in the theorem is the formal
dual to (and older than) the functor $\mathcal L$ in the preceding
section; it is due to Hausser and Nill \cite{HauNil:ITQHA}. We shall
recall some details in the proof.

\begin{Thm}\nmlabel{Theorem}{char}
  Let $H$ be a finite-dimensional quasibialgebra. The following
  are equivalent:
  \begin{enumerate}
    \item $H$ is a quasi-Hopf algebra.
    \item The monoidal category $\LModf H$ is rigid.
    \item The functor $\mathcal R\colon\LMod H\rightarrow\HMod{}HHH$
      is a category equivalence.
  \end{enumerate}
\end{Thm}
\begin{proof}
  The equivalence of (2) and (3) follows by duality from 
  \nmref{dualthm}. The implication (1)\imp(2) is the 
  motivation for the definition of a quasiantipode in Drinfeld's
  paper \cite{Dri:QHA}; the question whether one should require
  bijectivity of the antipode is actually irrelevant here.
  To recall some details: If $V$ is a finite-dimensional $H$-module,
  then a dual object for $V$ in $\LMod H$ is
  the dual vector space 
  $V^*$ with the $H$-module structure given by the transpose
  of the action via $S$, and the evaluation and coevaluation maps
  \begin{gather*}
    \ev\colon V^*\o V\ni\varphi\o v\mapsto \varphi(\alpha v)\in k\\
    \db\colon k\ni 1\mapsto \beta v_i\o v^i.
  \end{gather*}

  We shall now prove (3)\imp(1). We start by recalling the form
  of $\mathcal R$: it is the composition of 
  the functor $\LMod H\rightarrow\BiMod H$ which is given by
  restriction of the right module structure along
  $\epsilon\colon H\rightarrow k$, and the cofree right comodule
  functor $\BiMod H\rightarrow\HMod {}HHH$, where comodule now 
  means comodule over the coassociative coalgebra $H$ in the monoidal
  category $\BiMod H$. We observe that $\mathcal R$ is right adjoint,
  being the composition of two rather standard right adjoint functors:
  The left adjoint to the cofree comodule functor is just the 
  underlying functor $\HMod{}HHH\to\BiMod H$, while the left
  adjoint to the restriction functor $\LMod H\rightarrow\BiMod H$
  is the induction functor, in this case 
  mapping $M\in\BiMod H$ to $M/MH^+\in\LMod H$. 
  Thus we have
  the overall left adjoint $\mathcal F$ to $\mathcal R$ mapping
  $M\in\HMod {}HHH$ to $\ol M:=M/MH^+\in\LMod H$, and it is easy to find
  the unit of adjunction to be
  $$M\ni m\mapsto \ol{m\sw 0}\o m\sw 1\in M/MH^+\o H.$$
  In particular, 
  $$\ttheta\colon H_\str\o{_\str H^\str_\str} \ni g\o h\mapsto \ol{g\phi\so 1\o h\sw 1\phi\so 2}\o h\sw 2\phi\so 3\in \ol{H_\str\o {_\str H^\str_\str}}\o H$$
  is an isomorphism. While $\ttheta$ is a morphism 
  in $\HMod{}HHH$ by construction, where the structures are as
  indicated by
  the dots, we may observe that it is also a left $H$-module
  map with respect to another set of module structures: 
  the action of $H$ on the left tensor factor of $H\o H$ gives another
  $H$-bimodule ${_\str H_\str}\o H_\str$, and thus a left module
  $\ol{_\str H\o H}$; it is obvious that $\ttheta$ is a 
  module map with respect to these structures as well. In particular
  $\ol{_\str H\o H}^{\dim H}\cong H^{\dim H}$, so that we have
  $\ol{_\str H\o H}\cong H$ as left modules by Krull-Schmidt.
  Pick an isomorphism $\tgamma\colon\ol{H\o H}\rightarrow H$ of 
  left $H$-modules. Setting $\gamma(h)=\tgamma(\ol{1\o h})$
  we find $\tgamma(\ol{g\o h})=g\gamma(h)$. Next, we use the regular
  left $H\o H$-module structure on $H\o H$, which induces an $H\o H$-module
  structure on $\ol{H\o H}$; via $\tgamma$, we get an $H\o H$-module
  structure on $H$, such that the action of the left tensor factor is
  the regular module structure of $H$. In any such $H\o H$-module
  structure, the action of the right tensor factor has the form
  $h\circ g=gS(h)$ for some algebra antiendomorphism $S$ of $H$.
  Thus 
  $$g\gamma(\ell h)=\tgamma(\ol{g\o \ell h})=\tgamma(\ol{g\o h})S(\ell)=g\gamma(h)S(\ell)$$
  for all $g,h,\ell\in H$. In particular
  $\gamma(h)=\beta S(h)$ for $\beta:=\gamma(1)$.
  We define $\vartheta:=(\tgamma\o H)\ttheta\colon H\o H\rightarrow H\o H$
  and find
  $\vartheta(g\o h)=g\phi\so 1\beta S(h\sw 1\phi\so 2)\o h\sw 2\phi\so 3$.
  Note that 
  $$\vartheta\colon H_\str\o {_\str H_\str^\str}\rightarrow {_\str(_SH)}\o{_\str H^\str_\str}$$
  is a morphism in $\HMod {}HHH$ with the indicated structures,
  where $_SH$ denotes the left $H$-module structure on $H$ given by
  $S$. 
  In addition, $\vartheta$ is an $H$-module map with respect
  to the left $H$-module structures given by the regular action
  of $H$ on the left tensor factors. We may summarize the three
  variants of $H$-linearity in the formula
  $$\vartheta((g\o h)\xi(j\sw 1\o j\sw 2))=(g\o h\sw 2)\vartheta(\xi)(S(h\sw 1)\o j)$$
  for $g,h,j\in H$ and $\xi\in H\o H$, in which all multiplications
  are now in the algebra $H\o H$.

  As a first application
  $$h\sw 1\beta S(h\sw 2)=(H\o\epsilon)\vartheta(h\sw 1\o h\sw 2)
      =(H\o\epsilon)(\vartheta(1\o 1)(1\o h))=\beta\epsilon(h)$$
  for $h\in H$.
  Next, we set $\alpha:=(H\o\epsilon)\vartheta\inv(1\o 1)$, and find
  $$(H\o\epsilon)\vartheta\inv(g\o h)=(H\o\epsilon)((g\o 1)\vartheta\inv(1\o 1)(h\sw 1\o h\sw 2))=g\alpha h$$
  for all $g,h\in H$.
  This implies further
  $$S(h\sw 1)\alpha h\sw 2=(H\o\epsilon)\vartheta\inv(S(h\sw 1)\o h\sw 2)
     =(H\o\epsilon)((1\o h)\vartheta\inv(1\o 1))=\epsilon(h)\alpha,$$
  for $h\in H$, and
  $$1=(H\o\epsilon)\vartheta\inv\vartheta(1\o 1)
      =(H\o\epsilon)\vartheta\inv(\phi\so 1\beta S(\phi\so 2)\o\phi\so 3)
      =\phi\so 1\beta S(\phi\so 2)\alpha\phi\so 3.
  $$
  We can determine the inverse of 
  $\vartheta$ using that $H_\str\o{_\str H_\str^\str}$
  is the cofree right $H$-comodule within $\BiMod H$ over $H$,
  so that we have
  \begin{multline*}  
    \vartheta\inv(g\o h)=(H\o\epsilon)\vartheta\inv((g\o h)\sw 0)\o(g\o h)\sw 1
       \\=(H\o\epsilon)\vartheta\inv(gS(\phi\som 1)\o \phi\som 2h\sw 1)\o\phi\som 3h\sw 2
       \\=gS(\phi\som 1)\alpha\phi\som 2h\sw 1\o \phi\som 3h\sw 2.
  \end{multline*}       
  We find that
  \begin{multline*}
    1=(H\o\epsilon)\vartheta\vartheta\inv(1\o 1)
      =(H\o\epsilon)\vartheta(S(\phi\som 1)\alpha\phi\som 2\o \phi\som 3)
      \\=S(\phi\som 1)\alpha\phi\som 2\beta S(\phi\som 3),
  \end{multline*}
  which was the last axiom missing to show that $(S,\alpha,\beta)$
  is a quasiantipode.
  \end{proof}
\begin{Rem}\nmlabel{Remark}{strrem}
\begin{enumerate}
  \item The map $\vartheta$ occurs already in Drinfeld's paper
    \cite{Dri:QHA}. Its inverse is the proper quasi-Hopf analog
    of the canonical map $\kappa\colon H\o H\ni g\o h\mapsto gh\sw 1\o h\sw 2\in H\o H$
    for an ordinary bialgebra $H$. The canonical map $\kappa$ 
    is well-known to
    be a bijection if and only if $H$ has an antipode. Note, however,
    that both $\vartheta$ and $\vartheta\inv$ involve the quasiantipode,
    while bijectivity of the naively copied map $\kappa$ does not
    seem to have a relation to the question when a quasibialgebra
    $H$ is a quasi-Hopf algebra.
  \item As in \nmref{cqhacor} we see that
    $V_\str\o {_\str H^\str_\str}\cong {_SV}\o{_\str H^\str_\str}\in\HMod{}HHH$
    for every $V\in\RMod H$, when $H$ is a finite-dimensional
    quasi-Hopf algebra. The same can be shown explicitly if 
    $H$ is not finite-dimensional, see \cite{HauNil:ITQHA,Sch:HMDQHA}.
\end{enumerate}
\end{Rem}
\section{Bijectivity of the antipode}\nmlabel{Section}{sec:bij}
Bulacu and Caenepeel have proved that the antipode (mapping) of
a finite-dimensional quasi-Hopf algebra $H$ is always a bijection.
As in the ordinary Hopf case, they show this along with the 
existence of integrals, and the fact that $H$ is a Frobenius
algebra. The standard textbook \cite{Swe:HA,Abe:HA,Mon:HAAR}
proof does not work in the quasi-Hopf situation: It is 
based on finding a Hopf module structure on the dual $H^*$.
A structure of Hopf module in $\HMod{}HHH$ is indeed given 
in \cite{HauNil:ITQHA} to do integral theory, but only
{\em using the assumption\/} that the antipode is bijective.

In this section we will give a rather short proof for:
\begin{Thm}[Bulacu--Caenepeel \cite{BulCae:IDQHAA}]
Let $H$ be a finite-dimensional quasi-Hopf algebra with antipode
$(S,\alpha,\beta)$. Then $H$ is a Frobenius algebra, and $S$ 
is bijective.
\end{Thm}
\begin{proof}Since $H\o H$ is an $H\o H$-$H\o H$-bimodule, we get an
$H\o H$-$H\o H$-bimodule structure on 
$$\Hom_{H\o k-}(H\o H,H):=\{F\colon H\o H\rightarrow H|\forall h\in H,\xi\in H\o H\colon F((h\o 1)\xi)=hF(\xi)\}$$
by setting
$((g\o h)F(j\o\ell)
  )(\xi)=F((1\o\ell)
  \xi(g\o h))j
$
for all $g,h,j,\ell\in H$, $F\in\Hom_{H\o k-}(H\o H,H)$
 and $\xi\in H\o H$. Thus
$\vartheta$ induces an automorphism $\tttheta$ of
$\Hom_{H\o k-}(H\o H,H)$ satisfying
$$\tttheta((S(g\sw 1)\o h)F(j\o g\sw 2))=(h\sw 1\o h\sw 2)\tttheta(F)(j\o g).$$

We have a bijection
$T\colon H\o H^*\rightarrow\Hom_{H\o k-}(H\o H,H)$
given by
$T(g\o\varphi)(h\o j)=h\varphi(j)g$.
It is straightforward to check that $T$ is an $H\o H$-$H\o H$-bimodule
map.

Thus $\tttheta$ induces an automorphism of $H\o H^*$ that is
an isomorphism between the left $H$-action on the right
tensor factor, and the diagonal left $H$-action.
The latter has the structure of a Hopf module 
${^\str_\str H_\str}\o{_\str H^*}\in\HMod H{}HH,$
so by the structure theorem for Hopf modules (which applies in the
left-right switched version since $H^\bop$ is a quasi-Hopf algebra),
it is a free $H$-module, isomorphic to $H^{\dim H}$.
By the Krull-Schmidt Theorem,  
$H^*\cong H$ as left $H$-modules. Thus $H$ is a Frobenius
algebra. Now we consider once more the isomorphism 
$\tttheta$, and identify $H\o H^*$ with $H\o H$ as left modules.
We see that $H_\str\o{_\str H}\cong{_\str H_\str}\o{_\str H}$
as $H$-$H$-bimodules. But 
${_\str H_\str}\o{_\str H}\cong {_\str H_\str}\o H_S$ by
the left-right switched version of part (2) of \nmref{strrem},
and it follows that we have an isomorphism of right $H$-modules
$$H\cong k\ou H(H_\str\o{_\str H})\cong k\ou H({_\str H_\str}\o H_S)\cong H_S.$$
Thus $S$ is a bijection.
\end{proof}
Our short proof took advantage of the structure theorem for Hopf
modules as well as the isomorphism $\vartheta$. It may be worthwhile
to note that one does not really need the full generality of
the structure theorem, but can use more directly the information 
contained
in the map $\vartheta$:

For $V\in\RMod H$, we can define
$$\vartheta_V\colon V\o H\ni v\o h\mapsto v\phi\so 1\beta S(h\sw 1\phi\so 2)\o h\sw 2\phi\so 3\in V\o H$$
and 
$$\vartheta_V\inv\colon V\o H\ni v\o h\mapsto vS(\phi\som 1)\alpha\phi\som 2h\sw 1\o \phi\som 3h\sw 2\in V\o H.$$
Since $\vartheta_V$ and $\vartheta_V\inv$ are natural in $V\in\RMod H$
and mutually inverse isomorphisms for $V=H$, we see that they are 
mutually inverse
isomorphisms for any $V\in\RMod H$. In particular, we see that
$V_\str \o{_\str H_\str}\cong {_SV}\o {_\str H_\str}$
as $H$-bimodules, so that every right $H$-module of the form
$V_\str\o H_\str $ is free. Since $H^\bop$ is a quasi-Hopf algebra, every
$H$-$H$-bimodule of the form $_\str H_\str \o{_\str V}$ with 
$V\in\LMod H$ is 
is isomorphic to $_\str H_\str\o {V_S}$, hence free 
as a left $H$-module. No other cases of the structure theorem for
Hopf modules were used in our proof.

\end{document}